\documentclass{article}
\usepackage{amssymb,amsmath,amsthm,longtable}
\usepackage[dvips]{graphicx}
\newcommand{\T}{{\mathbb T}}
\newcommand{\D}{{\mathbb D}}
\newcommand{\C}{{\mathbb C}}
\newcommand{\E}{{\mathbb E}}

\newcommand{\R}{{\mathbb R}}

\def\br#1{\left(#1\right)}
\def\brb#1{\left[#1\right]}
\def\brs#1{\left\{#1\right\}}
\newcommand{\K}{K}
\renewcommand{\L}{{\mathcal L}}
\renewcommand{\cap}{{\mathrm cap\:}}
 
\theoremstyle{remark}
\newtheorem{rem}{Remark}
 
\theoremstyle{remark}
\newtheorem{notation}{Notation}
 
\theoremstyle{definition}
\newtheorem{definition}{Definition}
 
\theoremstyle{plain}
\newtheorem{theorem}{Theorem}[section]

\newtheorem{lemma}[theorem]{Lemma} 
\newtheorem{corollary}[theorem]{Corollary} 
\date{}
\title{Random conformal snowflakes}
\author{D. Beliaev \and S. Smirnov }
 
\begin{document}
\maketitle
\begin{abstract}
In many problems of classical analysis extremal configurations
appear to exhibit complicated fractal structure. This makes it much 
harder to describe extremals and to attack such problems.
Many of these problems are related to the {\em multifractal analysis}
of {\em harmonic measure}. 

We argue that, searching for extremals 
in such problems, one should work with random fractals rather than 
deterministic ones.
We introduce a new class of fractals {\em random conformal snowflakes}
and investigate its properties developing tools to estimate spectra
and showing that extremals can be found in this class. As an
application we significantly improve known
estimates from below on the extremal behaviour of harmonic measure, 
showing how to constuct a rather simple snowflake, which has a spectrum 
quite close to the conjectured extremal value.

\end{abstract}

\section{Introduction}

It became apparent during the last decade that extremal configurations 
in many important problems in classical complex analysis exhibit 
complicated fractal structure. This makes such problems more difficult 
to approach than similar ones, where 
extremal objects are smooth. 

As an example one can consider the coefficient problem for univalent functions.
Bieberbach formulated his famous conjecture arguing that the K\"oebe 
function, which maps a unit disc to a plane with a stright slit, is extremal. 
The Bieberbach conjecture was ultimately proved by de Branges in 1985
\cite{deBranges}, while the sharp growth asymptotics was obtained by
Littlewood \cite{Littlewood25} in 1925 by a much easier argument.

However, coefficient growth problem for bounded functions remains 
widely open, largely due to the fact that the extremals must be of 
fractal nature (cf \cite{CaJo}). This relates (see \cite{BeSmECM}) 
to a more general question of finding the
{\em universal multifractal spectrum} of {\em harmonic measure}
defined below, which includes many other problems, in particular
conjectures of Brennan, Carleson and Jones,
Kraetzer, Szeg\"o, and Littlewood.

In this paper we report on our search for extremal fractals.
We argue that one should 
study  random fractals instead of deterministic ones. 
We introduce a new class of random fractals, {\it random conformal snowflakes},
investigate its properties, and as a consequence significantly improve
known estimates from below for the multifractal spectra of harmonic measure.

\subsection{Multifractal analysis of harmonic measure}
 
It became clear recently that appropriate language for many problems in
geometric function theory is given by the {\em multifractal analysis} 
of {\em harmonic measure}. 
The concept of multifractal spectrum of a measure was introduced by
Mandelbrot in 1971  in \cite{Mandelbrot72, Mandelbrot74} in two papers devoted to the
distribution of energy in a turbulent flow.
We use the definitions that appeared
in 1986 in a seminal physics paper \cite{HJKPS} by
Halsey, Jensen, Kadanoff, Procaccia, Shraiman
who tried to understand and describe scaling laws of
physical measures on different fractals of physical nature
(strange attractors, stochastic fractals like DLA, etc.).

There are various notions of spectra and several ways to make a
rigorous definition. Two standard spectra are {\em packing} and 
{\em dimension} spectra. The packing spectrum of harmonic measure
$\omega$ in a domain $\Omega$ with a compact boundary is defined as
$$
\pi_{\Omega}(t)=
\sup\,\Big\{q:\ \forall\delta>0~\exists~\delta-\mathrm{ packing }~\{B\}
~{\mathrm with}~\sum \mathrm{diam}(B)^t\omega(B)^q\,\ge\,1\Big\}\ ,
$$
where $\delta$-packing is a 
collection of disjoint open sets whose diameters do not exceed $\delta$.

The {\em dimension spectrum} 
which is defined in terms of harmonic measure $\omega$ on the boundary
of $\Omega$ (in the case of simply connected domain $\Omega$ harmonic
measure is the image under the Riemann map $\phi$ of the normalised length on the unit circle).
Dimension spectrum gives the dimension of the set of points, where
harmonic measure satisfies a certain power law:
$$
f(\alpha)~:=~\mathrm{dim}\,
\Big\{z:~\omega\br{B(z,\delta)}\,\approx\,\delta^\alpha\,,~\delta\to0\Big\},~\alpha\ge\frac12~.
$$
Here $\mathrm{dim}$ stands for the Hausdorff or Minkowski dimension,
leading to possibly different spectra. The restriction $\alpha \ge 1/ 2$ is due
to Beurling's inequality. 
Of course in  general  there will be many points
where measure behaves differently at different scales, so
one has to add $\limsup$'s and $\liminf$'s to the definition above
-- consult \cite{Makarov} for details.

In our context it is more suitable to work with a modification of
the packing spectrum which is specific for the harmonic measure on a two dimensional
 simply connected domain $\Omega$. In this case we can define 
 the {\it integral means spectrum} as
$$
\beta_\phi(t)~:=~\limsup_{r\to1+}\frac{\log \int_{0}^{2\pi}|\phi'(re^{i\theta})|^t 
d\theta}{|\log(r-1)|},~t\in\R~,
$$
where $\phi$ is a Riemann map from the complement of the unit disc onto 
a simply connected domain $\Omega$.

Connections between all these spectra for particular domains
are not that simple, but the {\em universal spectra}
$$
\Pi(t)=\sup_\Omega \pi(t), \quad
F(\alpha)=\sup_\Omega f(\alpha), \ \ \mathrm{and} \ \ B(t)=\sup_\Omega\beta(t)
$$
are related by Legendre-type transforms:
\begin{eqnarray*}
F(\alpha)&=&\inf_{0\le t \le 2} (\alpha \Pi(t)+t), \quad \alpha \ge 1~, \\
\Pi(t)&=&\sup_{\alpha\ge 1} \left(\frac{F(\alpha)-t)}{\alpha}\right), \quad 0\le 2\le  2~, \\
\Pi(t)&=&B(t)-t+1~.
\end{eqnarray*}
See Makarov's survey \cite{Makarov} for details.

\subsection{Random fractals}
 
One of the main problems in the computation of the integral means spectrum 
(or other multifractal spectra) is the fact that 
the derivative of a Riemann map for a fractal domain depends on the argument in a very 
non regular way: $\phi'$ is a ``fractal'' object in itself. 
We propose to study random fractals to overcome this problem. 
For a random function $\phi$
 it is natural to consider the
{\em average integral means spectrum:}
\begin{eqnarray*}
\bar\beta(t)&=&\sup\brs{\beta: \int_1(r-1)^{\beta-1}\int_0^{2\pi}
\E\brb{|f'(r e^{i\theta})|^t}d \theta d r=\infty}
\\
&=&\inf\brs{\beta: \int_1(r-1)^{\beta-1}\int_0^{2\pi}
\E\brb{|f'(r e^{i\theta})|^t}d \theta d r<\infty}.
\end{eqnarray*}
The average spectrum does not have to be related to the spectra of a 
particular realization. We want to point out that even if 
$\phi$ has the same spectrum a.s. it does not guarantee that $\bar\beta(t)$
equal to the a.s. value of $\beta(t)$. Moreover, it can
happen that $\bar\beta$ is not a spectrum of
{\em any} particular domain. 
 
But one can see that $\bar\beta(t)$ is bounded by 
the universal spectrum $B(t)$. Indeed, suppose that there is a random 
$f$ with $\bar\beta(t)>B+\epsilon$, hence for any $r$ there are
particular realizations of $f$ with $\int |f'(z)|d\theta>(r-1)^{-B-\epsilon/2}$. 
Then by Makarov's fractal approximation \cite{Makarov} there is a (deterministic) function
$F$ such that $\beta_F(t)>B(t)$ which is impossible by the definition of $B(t)$.
  
For many classes of random fractals $\E|\phi'|^t$ (or its growth rate) does not depend
on the argument. This allows us to drop the integration with respect to the argument
and study the growth rate along any particular radius. Perhaps more importantly 
$\E |\phi'|$ is no longer a ``fractal'' function.
 
One can think that this is not a big advantage compared
to the usual integral means spectrum: instead of averaging over different arguments we 
average over different realizations of a fractal. But most fractals are results of
some kind of an iterative construction, which means that they are invariant under 
some (random) transformation. Thus $\E|\phi'|^t$ is a solution of some kind of equation. 
Solving this equation (or estimating its solutions) we can find $\bar\beta(t)$.
 
In this paper we want to show how one can employ these ideas. 
In the Section \ref{sec:def} we introduce a new class of random fractals 
that we call random conformal snowflakes. In the Section \ref{sec:spectrum} 
we show that $\bar\beta(t)$ for this class is related to the main 
eigenvalue of a particular integral operator. We also prove 
the fractal approximation for this class in the Section \ref{sec:approximation}.
In the Appendix \ref{sec:application} we give an example of a snowflake 
and prove that for this snowflake $\bar\beta(1)>0.23$. This significantly improves 
previously known estimate $B(1)>0.17$ due to Pommerenke \cite{Pommerenke75}.

\section{Conformal snowflake}
\label{sec:def}
 
The construction of our conformal snowflake is
similar to the construction in Pommerenke's paper \cite{Pommerenke67lms}.
The main difference is the introduction of the randomness. 

By $\Sigma'$ we denote a class of all univalent functions 
$\phi:\D_-\to \D_-$ such that $\phi(\infty)=\infty$ and $\phi'(\infty)\in \R$.
Let $\phi\in \Sigma'$ be a function with expansion at infinity $\phi(z)=b_1 z+\dots$, 
by $\cap\phi=\cap \Omega$ we denote the 
logarithmic capacity of $\phi$ which is equal to $\log|b_1|$.  
We will also use the so called {\em Koebe $n$-root transform} 
which is defined as
$$
(\K\phi)(z)=(\K_n\phi)(z)=\sqrt[n]{\phi(z^n)}.
$$ 
It is a well known fact that the Koebe transform is well
defined and $\K\phi\in\Sigma'$.
It is easy to check that Koebe transform divides capacity by $n$ 
and that the capacity of a composition is the sum of capacities.
 
First we define the deterministic snowflake. 
To construct a snowflake we need a  building block $\phi \in \Sigma'$ 
and an integer $k\ge 2$.  
Our snowflake will be  the result of the following iterative procedure:
we start with the building block and at $n$-th step we take a composition
of our function and the $k^n$-root transform of the rotated building block.
 
\begin{notation}
Let $\phi \in \Sigma'$ and $\theta\in [0,2\pi]$. By $\phi_\theta(z)$ we 
denote the map whose range is the rotation of that for $\phi$, 
namely  $e^{i\theta}\phi(z e^{-i\theta})$.
\end{notation}
 
\begin{definition}
Let $\phi\in \Sigma'$, $k\ge 2$ be an integer number, and $\{\theta_n\}$ be a sequence 
of numbers from $\T$.  Let $f_0(z)=\phi_{\theta_0}(z)$ and 
\begin{eqnarray*}
f_{n}(z)=f_{n-1}(K_{k^n}\phi_{\theta_n}(z))=
\phi_{\theta_0}(\phi_{\theta_1}^{1/k}(\dots\phi_{\theta_n}^{1/k}(z^{k^n})\dots).
\end{eqnarray*}
The conformal snowflake $f$  is the limit of $f_n$.
For simplicity $S=\C\setminus f(\D_-)$ and $g=f^{-1}$ are also called a snowflake. 
\end{definition}

Sometimes it is easier to work with a slightly different symmetric snowflake
$$
\bar f_{n}(z)=\phi_{\theta_1}^{1/k}(\dots\phi_{\theta_n}^{1/k}(z^{k^n})\dots)=
\Phi_1\circ \cdots \circ\Phi_n(z),
$$
where $\Phi_j=\K_{k^j}\phi_{\theta_j}$. 
There are two equivalent ways to construct the symmetric snowflake from the usual one. 
One is to take the Koebe transform  $\K f_n$, another is to start with $f_0(z)=z$. 
It is easy to see that
$f_n=\Phi_0\circ\cdots\circ\Phi_n.$

How this snowflake grows? 
This is easy to analyse looking at the evolution of $\bar f_n$. At every step
we add $k^n$ equidistributed (according to the harmonic measure) small copies of the building block.
But they are  not exact copies, they are distorted a little bit by a conformal mapping.
 
\begin{figure}
\centering{
\begin{tabular}{ll}
\includegraphics[width=4cm,height=4cm]{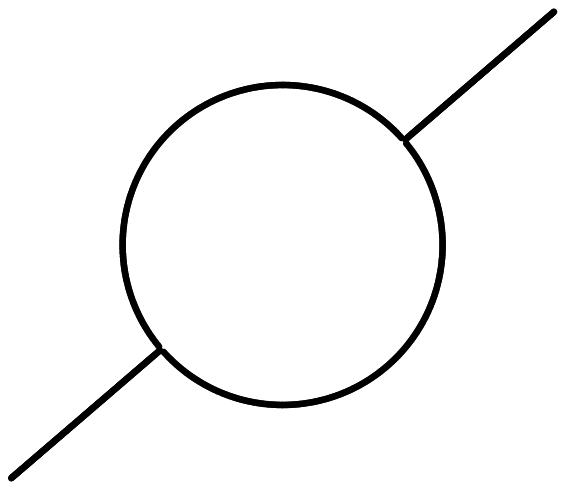} & 
\includegraphics[width=4cm,height=4cm]{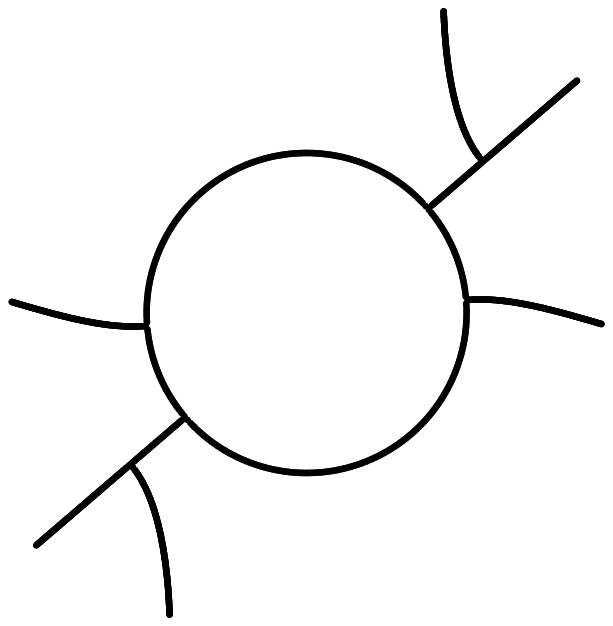} \\
\includegraphics[width=4cm,height=4cm]{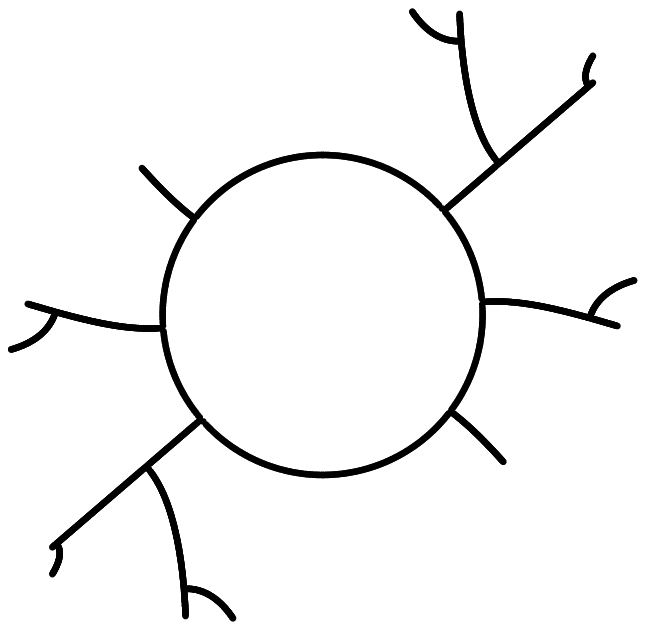} & 
\includegraphics[width=4cm,height=4cm]{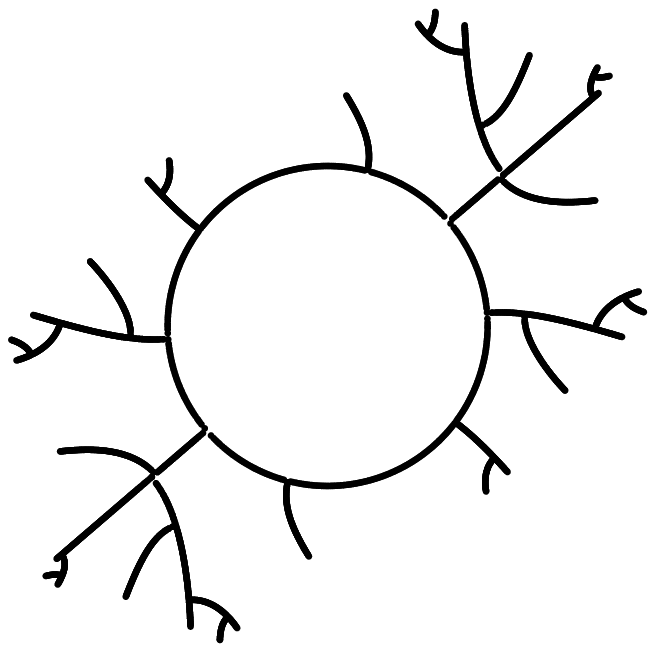} \\
\end{tabular}}
\caption{Images of $\bar f_1$, $\bar f_2$, $\bar f_3$, and $\bar f_4$}
\label{fbar}
\end{figure}

\begin{figure}
\centering{
\begin{tabular}{ll}
\includegraphics[width=4cm ]{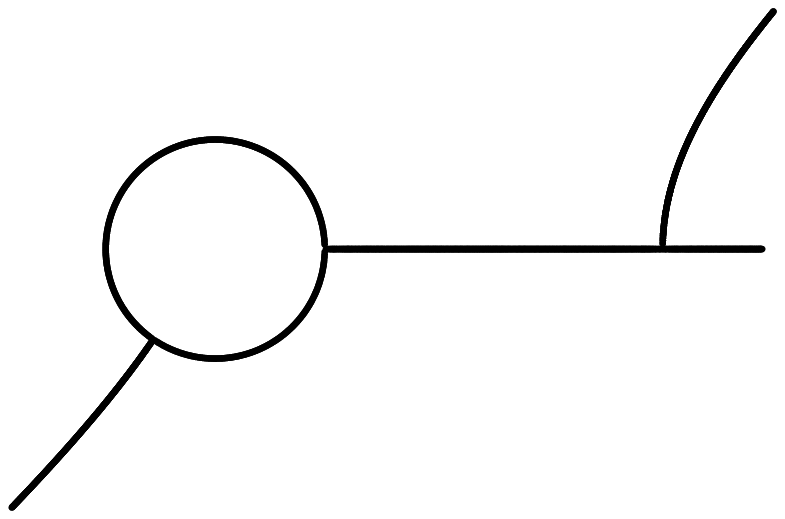} & 
\includegraphics[width=4cm ]{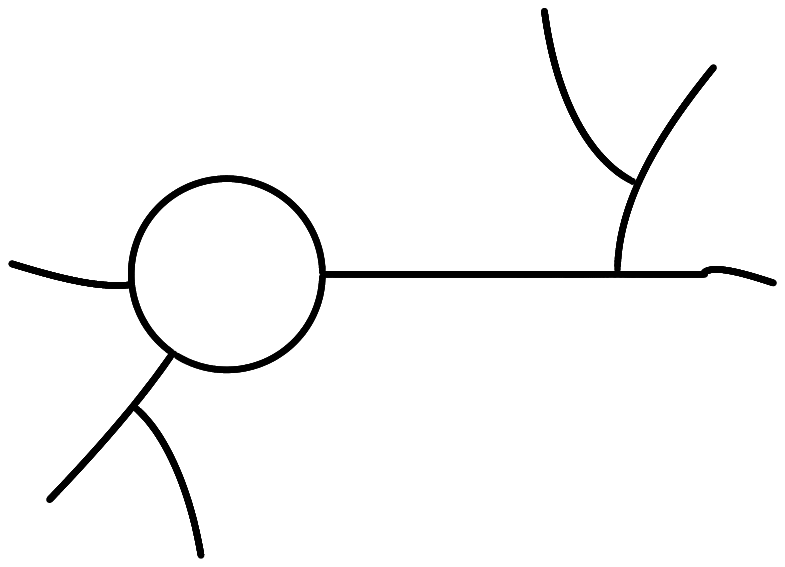} \\
\includegraphics[width=4cm ]{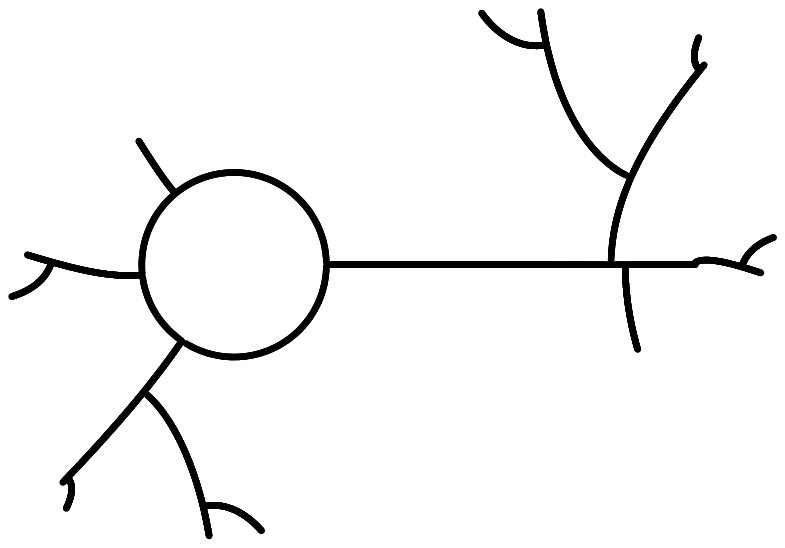} & 
\includegraphics[width=4cm ]{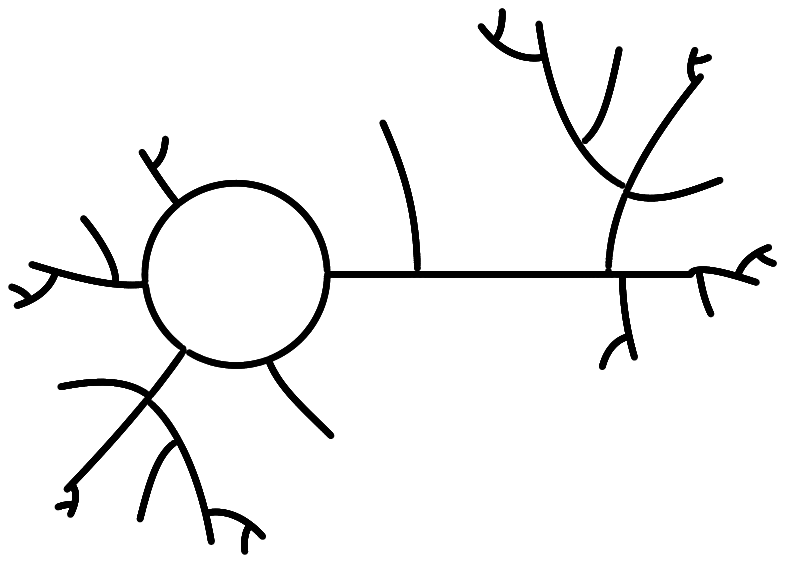} \\
\end{tabular}}
\caption{Images of $ f_1$, $  f_2$, $  f_3$, and $  f_4$}
\label{f}
\end{figure}

Figures \ref{fbar} and \ref{f} show images of the first four functions $\bar f$ and $f$ 
with $k=2$ and the building block is a slit map (which adds a straight slit   of length $4$). 
 
\begin{figure}
	\centering
		\includegraphics{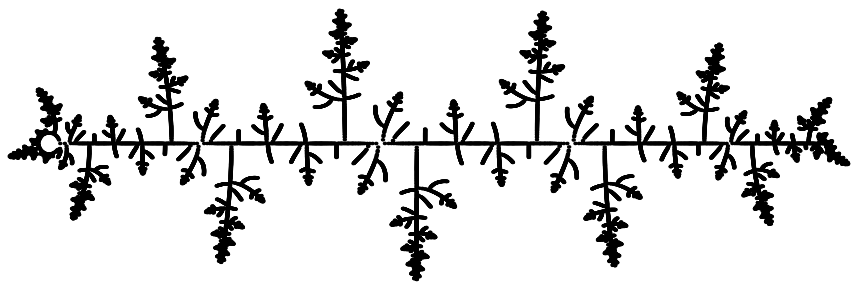}
		\caption{The third generation of a snowflake: $f_3$.}
	\label{pic2}
\end{figure}
 
\begin{figure}
	\centering
		\includegraphics{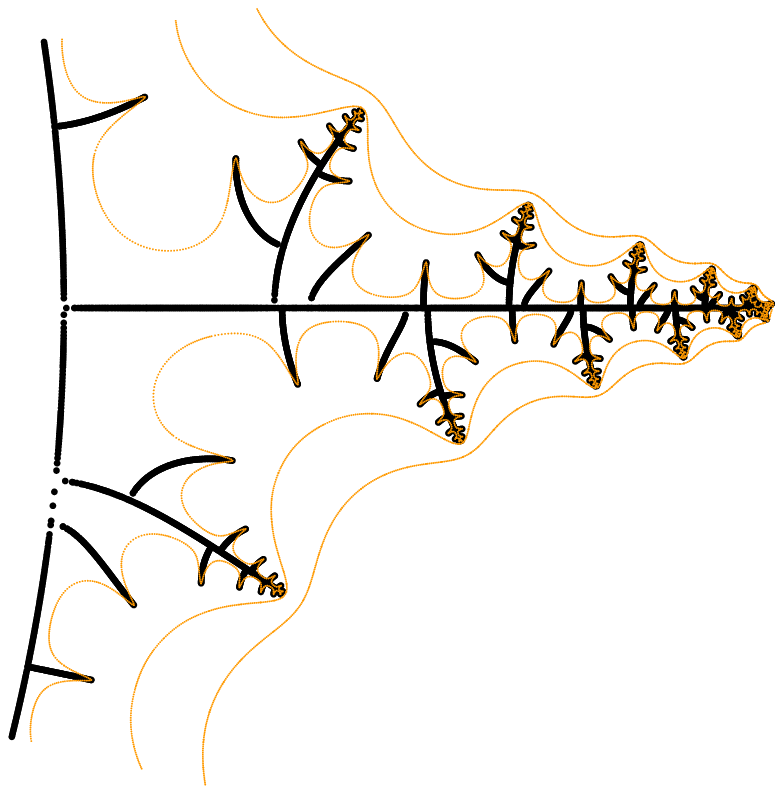}
	\caption{The image of a small boundary arc under $\bar f_3$ with three Green's lines.}
	\label{pic3}
\end{figure}

\begin{lemma} Let $f_n=\phi_{\theta_0}(\bar f_n(z))$ be the 
$n$-th approximation to the snowflake with a 
building block $\phi$ and $k\ge 2$. Then $\cap( f_n)$ and 
$\cap(\bar f_n)$ are bounded by (and converge to) $\cap(\phi)k/(k-1)$ 
and $\cap(\phi)/(k-1)$.
\label{l:boundedcapacity}
\end{lemma}
 
\begin{proof}
This lemma follows immediately from the standard facts that
\begin{eqnarray*}
\cap(f\circ g)&=&\cap(f)+\cap(g),\\
\cap(K_n f)&=&\cap(f)/n.
\end{eqnarray*}
\end{proof}
 
\begin{theorem} The conformal snowflake is well defined, namely
let $f_n$ be an $n$-th approximation to a snowflake with a building block $\phi$ and $k\ge 2$. 
Then there is $f\in\Sigma'$ such that $f_n$ converge to $f$
uniformly on every compact subset of $\D_-$.
\label{thm:conv}
\end{theorem}

\begin{proof}
Fix $\epsilon>0$. It is enough to prove that $\bar f_n$ converge uniformly on
$\D_\epsilon=\{z: |z|\ge 1+\epsilon\}$.
Suppose that $m>n$ so we can write 
$\bar f_m=\bar f_n \circ \Phi_{n,m}$ where 
$\Phi_{n,m}=\Phi_{n+1}\circ\cdots\circ \Phi_m$ and 
$$
|\bar f_n(z)-\bar f_m(z)|=
|\bar f_n(z)-\bar f_n(\Phi_{n,m}(z))|\le 
\max_{\zeta \in \D_\epsilon}|\bar f_n'(\zeta)||z-\Phi_{n,m}(z)|.
$$
By the Lemma \ref{l:boundedcapacity} $\cap(\bar f_n)$ is uniformly bounded, hence 
by the growth theorem the derivative of $\bar f_n$ is uniformly bounded in $\D_\epsilon$. 
Thus it is enough to prove that
$\Phi_{n,m}(z)$ converge uniformly to $z$. 
 
Let $\phi(z)=b_1 z+\dots$ at infinity, then 
$\Phi_n(z)=b_1^{1/k^{n}}z+\dots$ 
Functions
$\Phi_{n,m}$ have expansion 
$$
b_1^{k^{-n}+\dots +k^{-m}} z+\dots=b_1^{(n,m)}z+\dots
$$
 Obviously, $b_1^{(n,m)}\to 1$ as $n\to\infty$. 
This proves that $\Phi_{n,m}(z)\to z$ uniformly on $\D_\epsilon$ hence $f_n$
converge uniformly. Uniform limit of a functions from $\Sigma'$ 
can be either a constant or a function from $\Sigma'$. Since $\cap(f_n)$ is uniformly bounded 
the limit can not be a constant.
\end{proof}
 
\begin{definition}
Let $\phi\in\Sigma'$ and $k\ge 2$ be an integer number. The random conformal 
snowflake is a conformal snowflake defined by $\phi$, $k$, and $\{\theta_n\}$, where
$\theta_n$ are independent random variables uniformly distributed on $\T$.
\end{definition}
 
\begin{theorem}
\label{stationary} 
Let $\phi\in\Sigma'$,  $k\ge 2$ be an integer number, and $\psi=\phi^{-1}$. Let $f$ be a corresponding 
random snowflake and $g=f^{-1}$. Then the distribution of $f$ is invariant 
under the transformation $\Sigma'\times \T \mapsto \Sigma'$ which is defined by
$$
(f,\theta)\mapsto \phi_\theta(\K_k f). 
$$  
In other words
\begin{eqnarray}
f(z)&=&\phi_\theta((\K_k f)(z))= \phi_\theta(f^{1/k}(z^k)),\\
g(z)&=&(\K_k g)(\psi_\theta(z))= g^{1/k}(\psi^k_\theta(z)),
\end{eqnarray}
where $\theta$ is uniformly distributed on $\T$. 
Both equalities should be understood in the sense of distributions, i.e.
 distributions of both parts are the same.
\end{theorem}

\begin{proof}
Let $f$ be a snowflake generated by $\{\theta_n\}$. The probability distribution of 
the family of 
snowflakes is the infinite product of (normalised) Lebesgue measures on $\T$. 
By the definition
$$
f(z)=
\lim_{n\to\infty} \phi_{\theta_0}(\phi_{\theta_1}^{1/k}(\dots\phi_{\theta_n}^{1/k}(z^{k^n})\dots)
$$
and
$$
\phi_\theta((\K_k f)(z))=\lim_{n\to\infty} 
\phi_{\theta}(\phi_{\theta_0}^{1/k}(\dots\phi_{\theta_n}^{1/k}(z^{k^{n+1}})\dots),
$$
hence $\phi_\theta(\K_k f)$ is just the snowflake defined by the sequence 
$\theta, \theta_0, \theta_1,\dots$. So the transformation 
$f(z)\mapsto\phi_\theta((\K_k f)(z))$ is just a shift on the $[-\pi,\pi]^{\mathbb N}$. 
Obviously the product measure is invariant under the shift. This proves stationarity of $f$. 
Stationarity of $g$ follows immediately from stationarity of $f$.
\end{proof}

There is another way to think about random snowflakes. 
Let $\mathcal M$ be a space of probability measures on $\Sigma'$. 
And let $T$ be a random transformation $f \mapsto \phi_\theta (\K_k f)$, 
where $\theta$ is uniformly distributed on $[\pi,\pi]$. 
Obviously $T$ acts on $\mathcal M$. The distribution of a random snowflake 
is the only measure which is invariant under $T$. 
In some sense the random snowflake is an analog of a Julia set: 
it semi-conjugates $z^k$ and $\psi_\theta^k$.

\begin{center} 
\vspace{1cm}
\begin{tabular}{ccc}
& \vspace{-1.5cm} $\stackrel{f}{\longleftarrow}$ & \includegraphics[width=0.7cm]{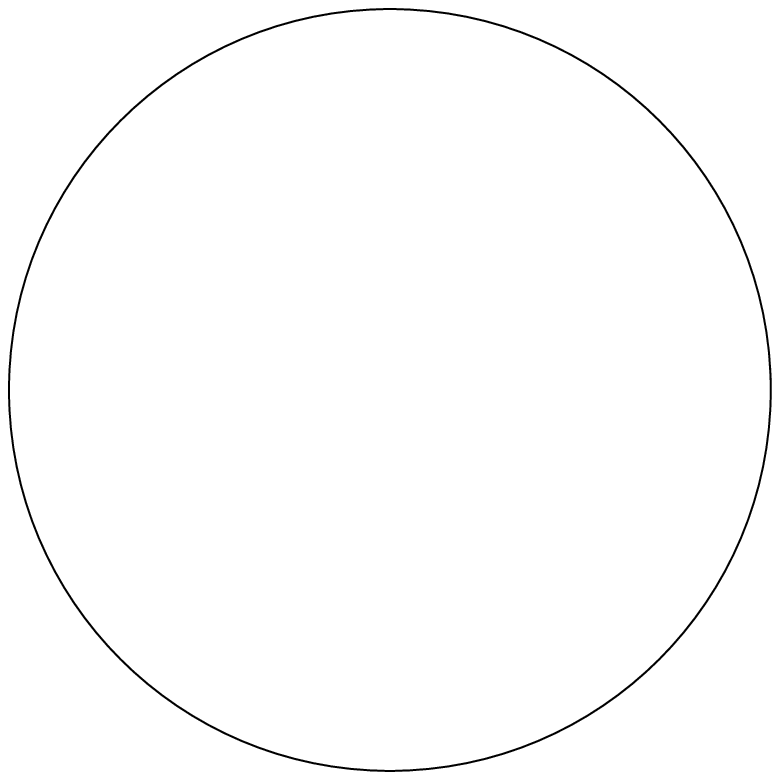} \\
\includegraphics[width=3cm]{snow24L.eps}&    &  \\
& & \\
$\uparrow \psi_\theta^k(z)$ & &   $\uparrow z^k$ \\
& \vspace{1cm} & \\
& \vspace{-1.5cm} $\stackrel{f}{\longleftarrow}$ & \includegraphics[width=0.7cm]{circle.eps}\\
\includegraphics[width=3cm]{snow24L.eps}&  & 
\end{tabular}
\end{center}

 The random conformal snowflakes are also rotationally invariant, the exact meaning 
is given by the following theorem.
 
\begin{theorem}
Let $\phi\in\Sigma'$, $k\ge 2$ and $g$ be the corresponding snowflake. Then $g$ is rotationally  
invariant, namely $g(z)$ and $e^{i\omega}g(e^{-i\omega}z)$ 
have the same distribution for any $\omega$.
\end{theorem}
 
\begin{proof}
Let $g_n(z)$ be the $n$th approximation to the snowflake defined by the sequence of 
rotations $\theta_0,\dots,\theta_n$. We claim that $\tilde g(z)=e^{i\omega}g(e^{-i\omega}z)$
is the approximation to the snowflake defined by $\tilde\theta_0,\dots,\tilde\theta_n$ where
$\tilde\theta_j=\theta_j+\omega k^j$ (we add arguments $\mod 2 \pi$). 
 
We prove this by induction. Obviously this is true for $\tilde g_0$. Suppose that 
it is true for $\tilde g_{n-1}$. By the definition of $g_n$ and assumption that
$g_{n-1}(e^{-i\omega}z)=e^{-i\omega}\tilde g_{n-1}(z)$ we have that
\begin{eqnarray*}
e^{i\omega}g_n(e^{-i\omega}z)&=&e^{i\omega}e^{i\theta_n/k^n}
\psi^{1/k^n}(e^{-i\theta_n} g_{n-1}^{k^n}(e^{-i\omega}z))
\\
&=&
e^{i\tilde\theta_n/k^n}\psi^{1/k^n}(e^{i\tilde\theta_n}\tilde g_{n-1}^{k^n}(z))
=
\tilde g_n(z).
\end{eqnarray*}
 
Obviously $\tilde \theta_n$ are also independent and uniformly 
distributed on $\T]$, hence $\tilde g_n$ has the same distribution as $g_n$.
\end{proof}
 
\begin{corollary}
The distributions of $|g(z)|$ and $|g'(z)|$ depend on $|z|$ only. The same is true for $f$.
\end{corollary}

\section{Spectrum of a conformal snowflake}
\label{sec:spectrum}
 
As we discussed above, for random fractals it is more natural to consider 
the average spectrum $\bar \beta(t)$ instead of
the usual spectrum $\beta(t)$. We will work with $\bar \beta(t)$ 
only and ``spectrum'' will always mean $\bar \beta(t)$.

\begin{notation}
 We will write $\L$ for the class of functions on $(1,\infty)$ that are bounded on 
compact sets and integrable in the neighbourhood of $1$. In particular, these functions belong to
$L^1[1,R]$ for any $1<R<\infty$. 
\end{notation}
 
Let $F(z)=F(|z|)=F(r)=\E \brb{\,|g'(r)/g(r)|^\tau\log^{\sigma}|g(r)|}$  where
$\tau=2-t$ and $\sigma=\beta-1$.
 
\begin{lemma} The $\bar\beta(t)$ spectrum of the snowflake is equal to
$$
\inf\brs{\beta: F(r)\in \L}.
$$
\end{lemma}
 
\begin{proof} 
By the definition $\bar\beta$ is the minimal value of $\beta$ such that
$$
\int_1\int_0^{2\pi}(r-1)^{\beta-1}\E|f(r e^{i\theta})|^t d\theta d r
$$ 
is finite. We change variable to $w=f(z)=f(re^{i\theta})$
\begin{eqnarray*}
\int\int \E\brb{|f'(re^{i\theta})|^t}(r-1)^{\beta-1} dr d\theta
&=&
\int \frac{\E\brb{|f'(z)|^t}(|z|-1)^{\beta-1}dm(z)}{r} 
\\&=&
\int \E\brb{\frac{|g'(w)|^{2-t}(|g(w)|-1)^{\beta-1}}{|g(w)|}}dm \, ,
\end{eqnarray*}
where $m$ is the Lebesgue measure. Note that $|g|$ is uniformly bounded and
$g$ is rotationally invariant, hence the last integral is finite 
if and only if 
$$
\int_1 |g'(r)|^\tau(|g(r)|-1)^\sigma d r<\infty.
$$
Since $1<|g|$ is uniformly bounded we have that
$|g'(r)|^\tau (|g(r)|-1)^\sigma$ is comparable up to an absolute constant to
$$
\br{\frac{|g'(r)|}{|g(r)|}}^\tau \log^\sigma|g(r)|.
$$
\end{proof}

\begin{lemma}If  $F \in \L$ then it
 is a solution of the following equation:
\begin{equation}
\label{eq}
F(r)=\frac{1}{k^\sigma}\int_{-\pi}^{\pi} F(|\psi^k(re^{i\theta})|)|\psi^{k-1}(re^{i\theta})
\psi'(re^{i\theta})|^\tau\frac{d\theta}{2\pi}.
\end{equation}
\end{lemma}
 
\begin{proof} 
By the Theorem \ref{stationary} $g(z)$ and $g^{1/k}(\psi_\theta^{k}(z))$ have the same distribution, hence
\begin{eqnarray*}
F(r)&=&\E\brb{\left|g'(r)/g(r)\right|^\tau\log^\sigma |g(r)|} \\
&=& \E \brb{\left|\frac{(g^{1/k}(\psi_\theta^k(r)))'}{g^{1/k}(\psi_\theta^k(r))}\right|^\tau 
\log^\sigma|g^{1/k}(\psi_\theta^k(r))|}\\
&=&
\E\brb{\left|\frac{g'(\psi^k_\theta(r))}{g(\psi_\theta^k(r))}\right|^\tau\log^\sigma(g(\psi^k_\theta(r)))
\frac{|\psi'_\theta(r)\psi^{k-1}_\theta(r)|^\tau}{k^\sigma}},
\end{eqnarray*}
where $\theta$ has a uniform distribution. 
The expectation is the integral with respect to the joint distribution
of $g$ and $\theta$, since they are independent this joint distribution
is just a product measure. So we can write it as the double integral: 
first we take the expectation with respect to the distribution of $g$ 
and than with respect to the (uniform) distribution of $\theta$
\begin{eqnarray*}
F(r)&=&\int_{-\pi}^\pi \br{\int  
\left|\frac{g'(\psi^k_\theta(r))}{g(\psi_\theta^k(r))}\right|^\tau\log^\sigma(g(\psi^k_\theta(r))) 
\frac{|\psi'_\theta(r)\psi^{k-1}_\theta(r)|^\tau}{k^\sigma}d g} \frac{d\theta}{2\pi}
\\
&=& \int_{-\pi}^\pi \br{\int  \left|\frac{g'(\psi^k_\theta(r))}
{g(\psi_\theta^k(r))}\right|^\tau\log^\sigma(g(\psi^k_\theta(r))) 
d g}\frac{|\psi'_\theta(r)\psi^{k-1}_\theta(r)|^\tau}{k^\sigma} \frac{d\theta}{2\pi}.
\end{eqnarray*}
 The inner integral is equal to $F(\psi_\theta^k(r))=F(\psi^k(e^{-i\theta}r))$ by the definition of $F$, hence 
\begin{eqnarray*}
F(r)&=&\int_{-\pi}^\pi F(\psi^k(e^{-i\theta}r))
\frac{|\psi'(e^{-i\theta}r)\psi^{k-1}(e^{-i\theta}r)|^\tau}{k^\sigma} \frac{d\theta}{2\pi}
\\
&=& \frac{1}{k^\sigma}\int_{-\pi}^\pi F(\psi^k(e^{i\theta}r))
|\psi'(e^{i\theta}r)\psi^{k-1}(e^{i\theta}r)|^\tau \frac{d\theta}{2\pi}
\end{eqnarray*}
which completes the proof. 
\end{proof}

This equation is the key ingredient in our calculations. One thinks about $F$ as the main 
eigenfunction of an integral operator. Hence the problem of finding the  
spectrum of the snowflake boils down to the question about the main eigenvalue
of a particular integral operator. Usually it is not very difficult to estimate the 
latter. 
 
 This justifies the definition:
\begin{equation}
Qf(r):=k \int_{-\pi}^{\pi} f (|\psi^k(re^{i\theta})|)\, |\psi^{k-1}(re^{i\theta})
\psi'(re^{i\theta})|^\tau\frac{d\theta}{2\pi}\, .
\end{equation}
Using this notation we can rewrite (\ref{eq}) as
$$
k^{\beta}F=QF.
$$
Note that this is in fact an ordinary kernel operator, $|\psi|$ is a smooth function 
of $\theta$, hence we can change the variable and write it as an integral operator. 
As mentioned above, the study of a $F$ is closely related to the study of operator $Q$
and its eigenvalues. And our estimate of the spectrum is in 
fact estimate of the main eigenvalue.
 
\subsection{Adjoint operator}
 
First of all we want to find a formally adjoint operator. 
Let $\nu$ be a bounded function and $R>1$ such that $D_R \subset \psi^k(D_R)$
where $D_R=\brs{z:1<|z|<R}$.
 
\begin{eqnarray*}
\int_1^R Qf(r)\nu(r)d r &=&
\int_1^R \nu(r)k \int_0^{2\pi} f(\psi^k(re^{i\theta}))\,
|\psi'(re^{i\theta})\psi^{k-1}(re^{i\theta})|^\tau\frac{d\theta}{2\pi}d r
\\
&=&
\int_{D_R}\frac{\nu(|z|)}{|z|}\frac{k}{ 2 \pi}f(\psi^k(z))
|\psi'(z)\psi^{k-1}(z)|^\tau d m(z),
\end{eqnarray*}
where $d m $ is the Lebesgue measure.
Changing a variable to $w=\psi^k(z)$ we get
\begin{eqnarray*}
&&\int_{\psi^k(D_R)}\frac{\nu(\phi(w^{1/k}))}{|\phi(w^{1/k})|}
\frac{1}{k 2\pi}f(w)|\phi'(w^{1/k})w^{1/k-1}|^{2-\tau} d m(w)
\\
&\ge&
\int_1^R\int_0^{2\pi k}
\frac{r\nu(\phi(r^{1/k}e^{i\theta/k}))}{|\phi(r^{1/k}e^{i\theta/k})|}
f(r)r^{\frac{(1-k)(\tau-2)}{k}}|\phi'(r^{1/k}e^{i\theta/k})|^{2-\tau}
\frac{d\theta}{2\pi k}d r
\\
&=&
\int_1^R f(r)\int_0^{2\pi }
\frac{r\nu(\phi(r^{1/k}e^{i\theta}))}{|\phi(r^{1/k}e^{i\theta})|}
r^{\frac{(1-k)(\tau-2)}{k}}|\phi'(r^{1/k}e^{i\theta})|^{2-\tau}
\frac{d\theta}{2\pi}d r.
\end{eqnarray*}
 
So we define another operator 
\begin{equation}
P\nu(r):={r^{1-\frac{(k-1)(2-\tau)}{k}}}
\int_0^{2\pi} \frac{\nu(\phi(r^{1/k}e^{i\theta}))}{|\phi(r^{1/k}e^{i\theta})|}
|\phi'(r^{1/k}e^{i\theta})|^{2-\tau}\frac{d\theta}{2\pi}.
\label{defP}
\end{equation}
 
Changing $2-\tau$ to $t$ we can rewrite  (\ref{defP}) as 
 
\begin{equation}
P\nu(r):={r^{1-\frac{(k-1)t}{k}}}
\int_0^{2\pi} \frac{\nu(\phi(r^{1/k}e^{i\theta}))}{|\phi(r^{1/k}e^{i\theta})|}
|\phi'(r^{1/k}e^{i\theta})|^{t}\frac{d\theta}{2\pi}.
\end{equation}
 
The inequality above can be written as
\begin{equation}
\label{cover}
\int_1^R Q f(r)\nu(r)d r \ge
\int_1^R f(r)P\nu(r) d r.
\end{equation}
We would like to note that for $R=\infty$ there is an equality since $\psi^K(\D_-)$ 
covers $\D_-$ exactly $k$ times. In this case 
\begin{equation}
\label{conj}
\int_1^\infty Q f(r)\nu(r)d r =
\int_1^\infty f(r)P\nu(r) d r.
\end{equation}
so operators $P$ and $Q$ are formally adjoint on $[1,\infty)$.

\begin{lemma}
Operator $Q=Q(t)$ acts on $\L$ if $\int|\phi'(r e^{i\theta})|d \theta$ is bounded. 
If $t\ge 1$ then it also acts on $L^1(1,\infty)$. 
\end{lemma}
 
\begin{proof}
Let $\nu=1$ in (\ref{conj}). Then
$$
\int_1^\infty Q f(r)d r =
\int_1^\infty f(r)
\int_{-\pi}^{\pi} \frac{r^{1-\frac{(k-1)t}{k}}}{|\phi(r^{1/k}e^{i\theta})|}
|\phi'(r^{1/k}e^{i\theta})|^{t}\frac{d\theta}{2\pi}.
$$
Let $r<R$, then 
$$
\frac{r^{1-\frac{(k-1)t}{k}}}{|\phi(r^{1/k}e^{i\theta})|}<R^{1-\frac{(k-1)t}{k}},
$$ 
so the second integral is bounded since   $\int |\phi'|^t d\theta$ is bounded. 
This proves that $Q f$  is in $\L$. To prove that it acts on $L^1(1,\infty)$ we should 
consider the  large values of $r$.
At infinity $\phi(z)=c z+ \dots $ and $\phi'(z)=c+\dots$, hence
$$
\begin{aligned}
\frac{r^{1-\frac{(k-1)t}{k}}}{|\phi(r^{1/k}e^{i\theta})|}
|\phi'(r^{1/k}e^{i\theta})|^{t}\approx \frac{r^{1-\frac{(k-1)t}{k}} |c|^t}{|c|r^{1/k}}=
\mathrm{const}\, r^{1-\frac{(k-1)t)}{k}-\frac{1}{k}}
\\
=\mathrm{const}\, r^{\frac{(k-1)(1-t)}{k}}.
\end{aligned}
$$
Thus the second integral is comparable (up to a universal constant) to $r^{\frac{(k-1)(1-t)}{k}}$, so it is bounded if and only if $t\ge 1$. 
\end{proof}
 
\begin{rem}
Note that the assumption on the integral of $|\phi'|$ is just a bit stronger than 
$\beta_\phi(t)=0$. We restrict ourselves to the building blocks that
are smooth up to the boundary, for such building blocks this assumption
is always true.
Condition $t\ge 1$ is technical and due to the behavior at infinity which should be irrelevant. 
Introducing the weight at infinity we can get 
rid of this assumption. 
\end{rem}
 
Next we want to discuss how eigenvalues of $P$ and $Q$ are related to 
the spectrum of the snowflake.
If $F$ is integrable then it is  a solution of (\ref{eq}) and using (\ref{cover}) 
we can write
\begin{equation}
\int_1^R F(r)\nu(r) =
\int_1^R \frac{Q F (r)}{k^{1+\sigma}}\nu(r)\ge
\int_1^R F(r) \nu(r) \frac{P\nu(r)}{\nu(r)k^{\sigma+1}}.
\end{equation}
Suppose that $t$ is fixed. Let us fix a positive test function $\nu$. If   
$P\nu(r)\ge \nu(r)k^{\sigma+1}$ then we arrive at contradiction, this means 
that  $F(r)$ for this particular pair of $\tau$ and $\sigma$ can not be integrable.
Using this fact we can estimate $\bar\beta(t)$ from below. Hence any positive 
$\nu$ gives the lower bound of the spectrum.
\begin{equation}
\label{betalog}
\bar\beta(t)\ge \min_{1\le r \le R}\log\left(\frac{P\nu(r)}{\nu(r)}\right)/\log k.
\end{equation}
Obviously, the best choice of $\nu$ is an eigenfunction of $P$ corresponding 
to the maximal eigenvalue. This proves the following lemma:
 
\begin{lemma}
Let $\lambda$ be  the maximal eigenvalue of $P$ (on any interval $[1,R]$ such that 
$D_R \subset \psi^k(D_R)$) then $\bar\beta(t)\ge \log \lambda/\log k$.
\end{lemma}
 
\section{Fractal approximation}
\label{sec:approximation}
 
In this section we prove the fractal approximation by
conformal snowflakes. Namely we show that for any $t$ one can construct a
snowflake with building block which is smooth up to the boundary
and $\bar\beta(t)$ arbitrary close to $B(t)$. The proof of this theorem is
similar to the proof of the fractal approximation for standard snowflakes but it
is less technical. 
 
\begin{theorem}
\label{thm:approximation}
For any $\epsilon$ and $t$
there are a building block $\phi\in\Sigma' \cap C^\infty(\{|z|\ge 1\})$ 
and a positive integer $k$ that define the snowflake with 
$\bar\beta(t)>B(t)-\epsilon$.
\end{theorem}
 
We will use the following lemma.
\begin{lemma}
\label{l:polygon}
For any $\epsilon>0$, $t\in \R$ there is $A>0$ such that for any  $\delta>0$ there is a  function 
$\phi\in\Sigma' \cap C^\infty $ such that
$$
\int \left|\phi'(re^{i\theta})\right|^t d \theta >
A \br{ \frac{1}{\delta}}^{B(t)-\epsilon}
$$
for $\delta>r-1$.  Moreover, capacity of $\phi$ is bounded by a universal constant 
that does not depend on $\delta$.
\end{lemma}
 
\begin{proof}
 There is a function $f$ with $\beta_f(t)>B(t)-\epsilon$. Hence there is a constant $A$ such that 
$$
 \int \left| f'(r e^{i\theta}) \right|^t d\theta > A (r-1)^{-B(t)+2\epsilon}.
$$
The only problem is that this function is not smooth up to the boundary. Set $\phi(z)=f(s z)$. Obviously, $\phi(z) \rightrightarrows f(z)$ as $s \to 1$. If we fix a scale $\delta$ then there is  $s$ sufficiently close to $1$ such that 
$\int|\phi'/\phi|^t d \theta > A \delta^{-B(t)+2\epsilon}/2$.  But for $r<1+\delta$ the integral can not be smaller by the  subharmonicity. 
\end{proof}

\begin{proof}[Proof of the Theorem \ref{thm:approximation}.] 
It is easy to see that 
$$
\cap(f_n)<\cap(f)=\cap(\phi)/(1-1/k)<2\cap(\phi),
$$ 
hence $\cap(f_n)$ and $|f_n(z)|$ for $|z|<2$ are bounded by the universal constants that depend
on capacity of $\phi$ only and do not depend on $k$.
It also follows that $|K_k f_n(z)|<1+c/k$ for $|z|<2$ and
$c$ depending on $\cap(\phi)$ only.
 
Let us fix $t$ and let $\phi$ be a function from the Lemma \ref{l:polygon} for $\delta=c/k$.
By $I(f,\delta)$ we denote 
$$
\int_{-\pi}^\pi \left| f'(re^{i\theta}) \right|^t d \theta,
$$ 
where $r=\exp(\delta)$.
 
The $k$-root transform changes integral means in a simple way:
$$
I(K_k f,\delta/k)=\int \left|\frac{f'(r^k e^{i k\theta})}{f^{(k-1)/k}(r^k e^{i k\theta})}\right|^t r^{t(k-1)}d\theta.
$$
As we mentioned before, the capacity of the snowflake is bounded by the universal constant, 
hence $|f|$ can be bounded by a universal constant. Thus
$$
I(K_k f,\delta/k)>\mathrm{const}\, I(f,\delta).
$$ 
 
The function $f_{n+1}$ is a composition of a (random) function $\phi_\theta$ with
$Kf_n$. 
The expectation of $I(f_{n+1},1/k^{n+1})$ conditioned on $f_n$ is
\begin{eqnarray*}
\E\brb{I(f_{n+1},1/k^{n+1})\mid f_n}&=&
\int \int |\phi_\theta'(K_k f_n(re^{i\xi}))|^t|(K_k f_n)'(re^{i\xi})|^t d\xi d \theta 
\\ 
&=&\int |(K_k f_n)'(re^{i\xi})|^t \int |\phi'(e^{-i\theta} K_k f_n(re^{i\xi}))|^t d \theta d \xi,
\end{eqnarray*}
where $r=\exp(1/k^{n+1})$.
We know that  $|K_k f_n(re^{i\xi})|-1<c/k$. By our choice of $\phi$ 
$$
\int |\phi'(|K_k f_n(re^{i\xi})|e^{-i\theta})|^t d \theta> A \br{\frac{k}{c}}^{B(t)-\epsilon}.
$$
So
\begin{eqnarray*}
\E \brb{I(f_{n+1},1/k^{n+1})} &>& A \br{\frac{k}{c}}^{B(t)-\epsilon}
\E \brb{I(K f_n,1/k^{n+1})}
\\
&>&A \br{\frac{k}{c}}^{B(t)-\epsilon} \mathrm{const}\, \E\brb{ I(f_n,1/k^n)}.
\end{eqnarray*}
Applying this inequality $n$ times we obtain
\begin{equation*}
\E \brb{I(f_n,1/k^n)}>\mathrm{ const}^n\,  \br{\frac{k}{c}}^{n(B(t)-\epsilon)}.
\end{equation*}
So
$$
\frac{\log\E\brb{ I(f_n,1/k^n)}}{n\log k}>B(t)-2\epsilon 
$$
for sufficiently large $k$. This completes the proof.
\end{proof}
 
\section{Appendix: example of an estimate}
\label{sec:application}
 
The main purpose of this section is to show that
using conformal snowflakes it is not very difficult to
find good estimates. Particularly it means that if one 
of the famous conjectures mentioned in the introduction 
is wrong,
then it should be possible to find a counterexample.

In this section we will give an example of a simple snowflake and
estimate its spectrum at $t=1$. We could do essentially the same computations
for other values of 
$t$, but $B(1)$ is of special interest because it 
is related to the coefficient problem and Littlewood conjecture (see \cite{BeSmECM}
for details).

As a building block we use a very simple 
 function: a straight slit map. 
We use the following scheme: first we define 
a building block and this gives us the operator $P$. By (\ref{betalog}) 
any positive function $\nu$ gives us an estimate on the spectrum.  
To choose $\nu$, we find the first eigenvector of discretized operator $P$
and approximate it by a rational function. We compute $P\nu$ using Euler's 
quadrature formula and estimate the error term. The minimum of $P\nu/\nu$ gives us 
the desired estimate of $\beta(1)$. For $t=1$ we give the rigorous estimate
of the error term in the computation of $P\nu/nu$,  for other values of 
$t$ we give approximate values (computed with less precision) without 
any estimates of the error terms.

\subsection{Single slit domain}
 
We use a straight slit functions. First we define the basic slit function
\begin{equation}
\label{slit}
\phi(z,l)=\phi_l(z)=\mu_2\br{\frac{\sqrt{\mu_1^2(z s)+l^2/(4k+4)}}{\sqrt{1+l^2/(4l+4)}}},
\end{equation} 
where $s$ is a constant close to $1$, 
$\mu_1$ and $\mu_2$ are the M\"obius transformation that maps $\D_-$ onto
the right half plane and its inverse:
\begin{eqnarray*}
\mu_1(z)=\frac{z-1}{z+1},\\
\mu_2(z)=\frac{z+1}{z-1}.
\end{eqnarray*}
We also need the inverse function
\begin{equation}
\psi(z,l)=\psi_l(z)=\phi(z,l)^{-1}.
\end{equation}
The function $\phi$ first maps $\D_-$ onto the right half-plane, than we cut off
  a straight horizontal slit starting at the origin and map it back. The image $\phi_l(\D_-)$ is
$\D_-$ with a horizontal slit starting from $1$. The length of the slit is $l$.
The derivative of a slit map has a singularities at points that are mapped
to $1$. But if we take $s>1$ then these singularities are not in $\D_-$. We set
$s=1.002$.

We study the snowflake generated by $\phi(z)=\phi_{73}(z)$ with $k=13$
(numbers $13$ and $73$ are found experimentally). Figure \ref{pic2}  show 
the image of the unit circle under  $f_3$. The Figure \ref{pic3} 
shows the image of a small arc under $\bar f_3$ and three Green's lines.

 First we have to find the critical
radius $R$ such that $D_R\subset \psi^k(D_R)$. By symmetry of $\phi$, the
critical radius is the only positive solution of 
$$
\psi^k(x)=x.
$$
This equation can not be solved explicitly, but we can solve it numerically 
(we don't care about error term since we can take any greater value of $R$). 
The approximate value of $R$ is $ 76.1568$. To be on the safe side we fix 
$R=76.2$. The disc takes just a small portion of
$\psi(D_R)$ which means that there is a huge overkill in the inequality (\ref{cover}).

By (\ref{betalog}) any positive function $\nu$ gives a lower bound of 
spectrum. And this estimate is sharp when $\nu$ is the main eigenfunction of $P$.
So we have to find an ``almost'' eigenfunction of $P$.

\subsection{Almost eigenfunction of operator P}
 
Even for such a simple building block we can not find the eigenfunction
explicitly. Instead we look for some sort of approximation. The first
idea is to substitute integral operator $P$ by its discretized version. Here we 
use the  most simple and  quite crude approximation. 
 
Choose sufficiently large $N$ and $M$. Let $r_n=1+(R-1)n/N$ and $\theta_m=
2\pi m /M$. Instead of $P$ we have an $N\times N$ matrix with elements
$$
P_{n,n'}=\sum_{m} r_n^{1-t(k-1)/k}\frac{|\phi'(r_n^{1/k}e^{i\theta_{m}})|^t}
{|\phi(r_n^{1/k}e^{i\theta_{m}})|M},
$$
where summation is over all indexes $m$ such that $r_{n'}$ is the nearest point to 
$|\phi(r_n^{1/k}e^{i\theta_m})|$.
This defines the discretized operator $P_N$. Let $\lambda_N$ and $V_N$ be the main
eigenvalue and the corresponding eigenvector. A priori, $\lambda_N$ should converge
to $k^{\beta(t)}$, but it is not easy to prove and not clear how to find the rate of
convergence. But this crude estimate gives us the fast test whether the pair 
$\phi$ and $k$ defines a snowflake with large spectrum or not (this is the way how we 
found $k=13$ and $l=73$). 
 
Instead of proving convergence of $\lambda_N$ and estimating the error term we will
study $V_N$ which is the discrete version of the eigenfunction. We 
approximate $V_N$ by a rational function of a relatively small degree (in our case $5$), 
or by  any other simple function. In our case we find the rational function by the linear 
least square fitting. In any way we get a nice and simple function $\nu$ which is
 supposed to be close to the eigenfunction of $P$.

We would like to note that procedure, by which we obtained $\nu$, 
is highly non rigorous, but 
that does not matter since as soon as we have some explicit 
function $\nu$ we can plug it into
$P$ and get the rigorous estimate of $\beta$.

\begin{figure*}[h]
\centering
\includegraphics[width=8cm]{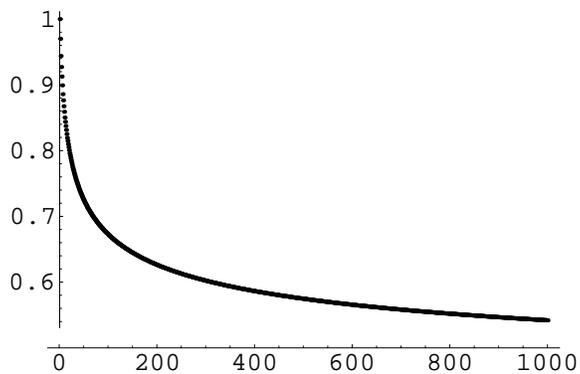}
\caption{Coordinates of the eigenvector}
\label{eigenvectors}
\end{figure*}

In our case we take $N=1000$ and $M=500$. The logarithm of the first 
eigenvalue is $0.2321$ (it is $0.23492$ if we take $s=1$). 
Figure \ref{eigenvectors} shows a plot with coordinates
of the first eigenvector.
 
\begin{figure}[h]
\centering
\includegraphics[width=8cm]{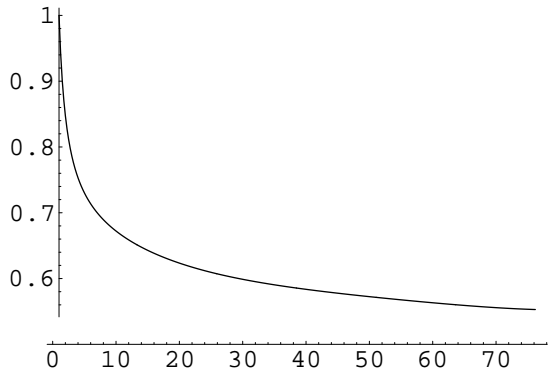}
\caption{An ``almost'' eigenfunction}
\label{nu}
\end{figure}

We scale this data from $[1,1000]$ to the interval $[1,R]$ and approximate by a rational function 
$\nu$:
$$
\begin{aligned}
\nu(x)=(7.1479+8.9280 x - 0.07765 x^2+ 1.733 \times 10^{-3} x^3 - 
\\ 2.0598 \times 10^{-5} x^4 + 
9.5353 \times 10^{-8}x^5)/( 2.7154+ 13.2845  x).
\end{aligned}
$$
 Figure \ref{nu} shows the plot of $\nu$. 
 
\subsection{Estimates of derivatives}

To estimate $\beta(1)$
we have to integrate $\nu(|\phi|)|\phi'|/|\phi|$. It is easy to see
that the main contribution to the derivative is given by a factor $|\phi'|/|\phi|$. 
Assume for a while that $s=1$. The fraction
$|\phi'(z)/\phi(z)|$ can be written as 
\begin{equation}
\label{fraction}
\frac{|z-1|}{|z|\sqrt{|(z-z_1)(z-z_2)}},
\end{equation}
 
where
$$
\begin{aligned}
z_1=\frac{-5033-292 i \sqrt{74}}{5625}\approx -0.894756 - 0.446556 i,
\\
z_2=\frac{-5033+292 i \sqrt{74}}{5625}\approx -0.894756 + 0.446556 i.
\end{aligned}
$$
 
Singular points $z_1$ and $z_2$ are mapped to $1$ and $\phi'$ has a square root type 
singularity at these points. They will play essential role in all further calculations.
We introduce notation $z_1=x+i y$ and $z_2=x-i y$, 
for $z$ we will use polar coordinates $z=r e^{i\theta}$.

We compute the integral of $f=|\nu(\phi)\phi'/\phi|$ using 
the Euler quadrature formula based on the trapezoid quadrature  formula 
$$
\int_0^{2\pi} f(x)dx \approx S_\epsilon^n(f)=
S_\epsilon(f) -\sum_{k=1}^{n-1} \gamma_{2k}\epsilon^{2k}\br{f^{(2k-1)}(2\pi)-f^{(2k-1)}(0)},
$$
where $S_\epsilon(f)$ is a trapezoid  quadrature formula with step $\epsilon$ and
$\gamma_k=B_k/k!$ were $B_k$ is the Bernoulli number.
The error term in the Euler formula is 
\begin{equation}
-\gamma_{2n}\max f^{(2n)}\epsilon^{2n}2\pi.
\label{error1}
\end{equation}
In our case function $f$ is periodic and terms with higher derivatives vanish. This means
that we can use (\ref{error1}) for any $n$ as an estimate of the error 
in the trapezoid quadrature formula.
 
Function $\phi$ has two singular points: $z_1$ and $z_2$. Derivative of $\phi$ blows up near 
these points. This is why we  introduce scaling factor $s$. We can write a power series of
$\phi$ near $z_1$ (near $z_2$ situation is the same by the symmetry)
$$
\phi^{(k)}=c_{-k}(z-z_1)^{-k+1/2}+c_{-k+1}(z-z_1)^{-k+3/2}+\dots+c_0+\dots\ .
$$
This means that for $s>1$ derivative can be estimated by
$$
|c_{-k}|(s-1)^{-k+1/2}+|c_{k+1}|(s-1)^{-k+3/2}+\dots\ .
$$
Tail of this series can be estimated because series converges in a disc of
a fixed radius (radius is $|z_1+1|$), and sum of tail can be estimated by a sum
of a geometric progression. Writing these power series explicitly we find
(for $s=1.002$)
\begin{alignat*}{3}
|\phi'| & <  55, & \quad  
|\phi''| & <  11800, & \quad
|\phi^{(3)}| &<  8.69\times 10^6,\\
|\phi^{(4)}| & <  1.08 \times 10^{10},  &\quad
|\phi^{(5)}| & <  1.90 \times 10^{13}, &\quad 
|\phi^{(6)}| & <  4.25 \times 10^{16}, \\
|\phi^{(7)}| & <  1.17 \times 10^{20}.  &\quad &  & &
\end{alignat*}

The maximal values for first six derivatives of $\nu$ are
\begin{alignat*}{3}
|\nu'|&<0.28, &\quad
|\nu''|&<0.45, &\quad
|\nu^{(3)}|&<1.12, \\
|\nu^{(4)}|&<3.69, &\quad
|\nu^{(5)}|&<15.3,  &\quad
|\nu^{(6)}|&<76.2.
\end{alignat*}
 
The derivative $\partial_\theta|\phi|$ can be estimated by $r |\phi'|$. We can write
sixth derivative of $\nu(|\phi|)|\phi'|/|\phi|$ as a rational function of 
partial derivatives of $|\phi|$, $|\phi'|$, and $\nu$.
 Than we apply triangle inequality and plug in the above estimates.
Finally we have
$$
\left|\frac{\partial\br{\frac{\nu(|\phi|)|\phi'|}{|\phi|}}}{\partial \theta^6}\right|<1.65\times 10^{21}.
$$
Plugging  $\epsilon=\pi/5000$ and estimate on sixth derivative into (\ref{error1})
we find that error term in this case is less than $0.0034$.
 
Next we have to estimate modulus of continuity with respect to $r$.
First we calculate 
$$
\partial_r|z-(a+b i)|^2=2 r - 2 (a\cos\theta+b\sin\theta).
$$
Applying this formula several times we find
$$
\begin{aligned}
\partial_r\br{\frac{|\phi'|}{|\phi|}}=
\partial_r \br{\frac{|z-1|}{r\sqrt{|z-z_1||z-z_2|}}} \le
\partial_r \br{\frac{|z-1|}{\sqrt{|z-z_1||z-z_2|}}} 
\\
=
\frac{r-\cos\theta}{|z-1|S}-\frac{r-x\cos\theta-y\sin\theta}{2|z-z_1|^2S}|z-1|-
\frac{r-x\cos\theta+y\sin\theta}{2|z-z_2|^2S}|z-1|,
\end{aligned}
$$
where $S=\sqrt{|z-z_1||z-z_2|}$.
Factoring out 
$$
\frac{1}{2|z-1|\cdot|z-z_1|^{5/2}|z-z_2|^{5/2}}
$$ we get
 
$$
\begin{aligned}
2(r-\cos\theta)|z-z_1|^2|z-z_2|^2
-(r-x\cos\theta-y\sin\theta)|z-1|^2|z-z_2|^2
\\
-(r-x\cos\theta+y\sin\theta)|z-1|^2|z-z_1|^2
\\
=-2(r^2-1)(2\cos^2\theta r (x-1)+\cos\theta (r^2+1)(x-1)+2 r y^2).
\end{aligned}
$$
This is a quadratic function with respect to $\cos\theta$. Taking values of $x$ and $y$ into
account we can write it as
$$
\cos^2\theta+\cos\theta\br{r+\frac{1}{r}}\frac{1}{2}-\frac{592}{5625}~.
$$
This quadratic function has two real roots. Their average is $-(r+1/r)/2<-1$, hence one
root is definitely less than $-1$. The product of roots is a small negative number, which meant that
the second root is  positive and less than $1$. Simple calculation shows that
this root decreases as $r$ grows. This means that the corresponding value of $\theta$
increases. Hence it attains its maximal value at $r=1.4$ and the maximal value is
at most $1.48$. This gives us that the radial derivative of $|\phi'|/|\phi|$ can be positive only
on the arc $\theta\in[-1.48,1.48]$. By subharmonicity it attains the maximal on the boundary of
$\{z\mid 1<r<1.4,\ -1.48<\theta<1.48\}$. It is not very difficult to check that 
maximum is at $z=1.4$ and it is equal to $0.36$. 
 
 Let 
$$
I(r)=\int_{-\pi}^\pi \nu(|\phi(r^{1/k}e^{i\theta})|
\left|\frac{\phi'(r^{1/k}e^{i\theta})}{\phi(r^{1/k}e^{i\theta})}\right| \frac{d \theta}{2\pi}.
$$
The derivative is
$$
I'(r)=\frac{1}{k r^{1-1/k}}\br{\int_{-\pi}^\pi
 \nu'(|\phi|)\partial_r |\phi| \frac{|\phi'|}{|\phi|}\frac{d \theta}{2\pi}+
\int_{-\pi}^\pi \nu(|\phi|)\partial_r\br{\frac{|\phi'|}{|\phi|}}\frac{d \theta}{2\pi}}.
$$
By the symmetry the first integral is zero. In the second integral 
$$
\nu(|\phi|)\partial_r\br{\frac{|\phi'|}{|\phi|}}
$$
can be positive only when $\theta\in[-1.48,1.48]$ and even in this case it is 
bounded by $0.36/(r^{1-1/k} k 2\pi)$. Hence
$$
I'(r)<2 \cdot 1.48 \cdot 0.36/(r^{1-1/k}k 2\pi)<0.0131 r^{1/k-1}.
$$
If we compute values $I(r_1)$ and $I(r_2)$ (with precision $0.0034$)
then the minimum of $P(\nu)/\nu$ on $[r_1,r_2]$ is at least  
\begin{equation}
r_1^{1/k}(\min\{I(r_1),I(r_2)\}-0.0034-0.0131(r_2-r_1) r^{1/k-1})/\nu(r_1).
\label{int}
\end{equation}
\begin{figure}[h]
\centering
\includegraphics{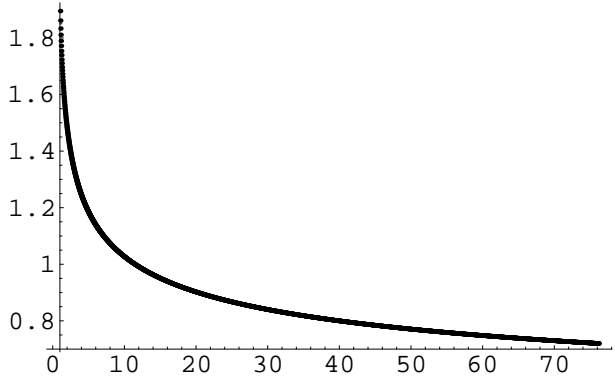}
\caption{Values of $I(r)$}
\label{integral}
\end{figure}
\begin{figure}[h]
\centering
\includegraphics{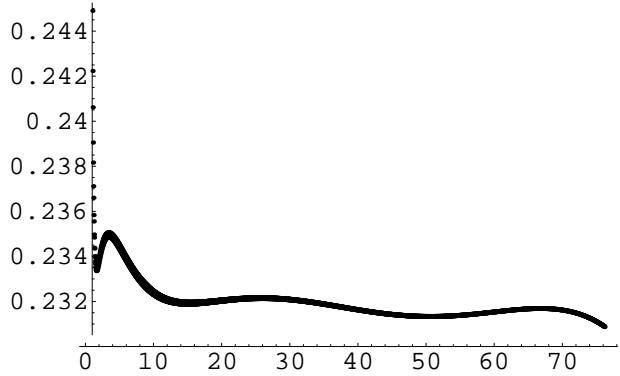}
\caption{Plot of $\log(P\nu/\nu)/\log k$}
\label{log}
\end{figure}
We take $3000$ equidistributed points on $[1,R]$ and compute $I(r)$ at these points.
The data for $I(r)$ is shown on the Figure \ref{integral}.
Applying the error estimate (\ref{int}) we find a rigorous estimate from below 
of $P\nu/\nu$. The minimum of $P\nu/\nu$ is at least $1.8079$ which means that
$$
\beta(1)>0.2308.
$$
The Figure \ref{log} shows the plot of $\log(P\nu/\nu)/\log k$.

\subsection{Estimates of spectrum for other values of $t$}
 
We also computed lower bound on spectrum of the same snowflake for other values of
$t$. Below are given base $13$ logarithms of eigenvalues of discretized operator $P$ ($N=1000$, $M=400$),
lower bounds of $\log(P\nu/\nu)/\log k$, $t^2/4$ and upper bound of the 
universal spectrum from \cite{HeSh, MaPo}. For values of $t$ close to zero we 
can not find a function that gives us a positive lower bound.

\medskip 

\begin{center}
\begin{tabular}{|l|l|l|l|l|} \hline 
\quad $t$ \qquad &\quad $\log_{13} \lambda$  
\qquad & \quad $\beta(t)>$ \qquad &\quad $t^2/4$\qquad & \quad $\beta(t)<$ \qquad\\ \hline
-2.0 & 0.6350 & 0.56 &1 & 1.218\\ \hline
-1.8& 0.5348& 0.48& 0.81&1.042 \\ \hline
-1.6&0.4395 &0.39 & 0.64& 0.871\\ \hline
-1.4&0.3502 &0.31 & 0.49& 0.706 \\ \hline
-1.2&0.2678 &0.220 & 0.36& 0.549\\ \hline
-1.0& 0.1936 &0.152 & 0.25& 0.403\\ \hline
-0.8&0.1290 &0.0925 & 0.16& 0.272\\ \hline
-0.6&0.0756 &0.0430 & 0.09& 0.159\\ \hline
-0.4& 0.0353& 0.0050& 0.04& 0.072\\ \hline
-0.2& 0.0100&0 & 0.01& 0.0179\\ \hline
0.2&0.0105 &0 & 0.01& 0.031\\ \hline
0.4&0.0387 &0.0280 & 0.04& 0.184\\ \hline
0.6& 0.0858& 0.0795& 0.09& 0.276\\ \hline
0.8& 0.1515&0.1505 & 0.16& 0.368\\ \hline
1.0&0.234 &0.234 & 0.25& 0.460\\ \hline
1.2&0.334 &0.332 & 0.36& 0.613\\ \hline
1.4& 0.448& 0.442& 0.49& 0.765\\ \hline
1.6& 0.576&0.570 &0.64 & 0.843\\ \hline
1.8&0.713 &0.698 & 0.81& 0.921\\ \hline
2.0& 0.859& 0.821& 1& 1\\ \hline
\end{tabular}
\end{center}
 
\bibliography{snow} 
\bibliographystyle{abbrv}
 
\end{document}